\magnification 1200
\advance\hoffset by 1,5truecm
\advance\hsize by 2,3 truecm
\def\makefootline{\baselineskip=52pt\line{\the\footline}}
\vsize= 23 true cm
\hsize= 14 true cm
\overfullrule=0mm

\font\grandsy=cmsy10 scaled \magstep0
\def\SS{{\grandsy x}}

\headline={\hfill\tenrm\folio\hfil}
\footline={\hfill}\pageno=1
\def\decale#1{\par\noindent\hskip 3em\llap{#1\enspace} \ignorespaces}
\def\undertext#1{$\underline{\hbox{#1}}$}

\newcount\coefftaille \newdimen\taille
\newdimen\htstrut \newdimen\wdstrut
\newdimen\ts \newdimen\tss

\def\fspeciale{\textfont0=\tenrmp%
\scriptfont0=\sevenrmp%
\scriptscriptfont0=\fivermp%
\textfont1=\tenip%
\scriptfont1=\sevenip%
\scriptscriptfont1=\fiveip%
\textfont2=\tensyp%
\scriptfont2=\sevensyp%
\scriptscriptfont2=\fivesyp%
\textfont3=\tenexp%
\scriptfont3=\tenexp%
\scriptscriptfont3=\tenexp%
\textfont\itfam=\tenitp%
\textfont\bffam=\tenbfp%
\textfont\slfam=\tenbfp%
\def\it{\fam\itfam\tenitp}%
\def\bf{\fam\bffam\tenbfp}%
\def\rm{\fam0\tenrmp}%
\def\sl{\fam\slfam\tenslp}%
\normalbaselineskip=12pt%
\multiply \normalbaselineskip by \coefftaille%
\divide \normalbaselineskip by 1000%
\normalbaselines%
\abovedisplayskip=10pt plus 2pt minus 7pt%
\multiply \abovedisplayskip by \coefftaille%
\divide \abovedisplayskip by 1000%
\belowdisplayskip=7pt plus 3pt minus 4pt%
\multiply \belowdisplayskip by \coefftaille%
\divide \belowdisplayskip by 1000%
\setbox\strutbox=\hbox{\vrule height\htstrut depth\wdstrut width 0pt}%
\rm}

\def\vmid#1{\mid\!#1\!\mid}

\def\fle{\rightarrow}

\null\vskip-1cm

\font\sc=cmcsc10
\font\grandsy=cmsy10 scaled\magstep2
\def\SSg{{\grandsy x}}
\newdimen\emm 
\def\pmb#1{\emm=0.03em\leavevmode\setbox0=\hbox{#1}
\kern0.901\emm\raise0.434\emm\copy0\kern-\wd0
\kern-0.678\emm\raise0.975\emm\copy0\kern-\wd0
\kern-0.846\emm\raise0.782\emm\copy0\kern-\wd0
\kern-0.377\emm\raise-0.000\emm\copy0\kern-\wd0
\kern0.377\emm\raise-0.782\emm\copy0\kern-\wd0
\kern0.846\emm\raise-0.975\emm\copy0\kern-\wd0
\kern0.678\emm\raise-0.434\emm\copy0\kern-\wd0
\kern\wd0\kern-0.901\emm}

\font\tendb=msbm10
\font\sevendb=msbm7

\newfam\dbfam
\textfont\dbfam=\tendb\scriptfont\dbfam=\sevendb\scriptscriptfont\dbfam=\sevendb
\def\db{\fam\dbfam\tendb}

\def\C{{\db C }}

\def\N{{\db N }}

\def\R{{\db R }}

\def\Z{{\db Z }}

\def\ie{\hbox{\it i.e.}}

\def\cotg{\mathop{\rm cotg}\nolimits}
\def\cos{\mathop{\rm cos}\nolimits}
\def\cotan{\mathop{\rm cotan}\nolimits}

%\font\bighelvetica=HelveticaB  at 14pt
%\font\abighelvetica=HelveticaB at 10pt
%\font\bbighelvetica=HelveticaB at 120pt
%\font\cbighelvetica=HelveticaB at 16pt
%\font\dbighelvetica=HelveticaB at 24pt
%\font\ebighelvetica=HelveticaB at 30pt

\font\gothique=eufm10
\def\got#1{{\gothique #1}}

%Pour une bibliographie
\newdimen\margeg \margeg=0pt
\def\bb#1&#2&#3&#4&#5&{\par{\parindent=0pt
    \advance\margeg by 1.1truecm\leftskip=\margeg
    {\everypar{\leftskip=\margeg}\smallbreak\noindent
    \hbox to 0pt{\hss\bf [#1]~~}{\bf #2 - }#3~; {\it #4.}\par\medskip
    #5 }
\medskip}}

\newdimen\margeg \margeg=0pt
\def\bbaa#1&#2&#3&#4&#5&{\par{\parindent=0pt
    \advance\margeg by 1.1truecm\leftskip=\margeg
    {\everypar{\leftskip=\margeg}\smallbreak\noindent
    \hbox to 0pt{\hss [#1]~~}{\pmb{\sc #2} - }#3~; {\it #4.}\par\medskip
    #5 }
\medskip}}

%Autre Modele avec le nom de l'auteur en sc mais non gras
\newdimen\margeg \margeg=0pt
\def\bba#1&#2&#3&#4&#5&{\par{\parindent=0pt
    \advance\margeg by 1.1truecm\leftskip=\margeg
    {\everypar{\leftskip=\margeg}\smallbreak\noindent
    \hbox to 0pt{\hss [#1]~~}{{\sc #2} - }#3~; {\it #4.}\par\medskip
    #5 }
\medskip}}

\def\messages#1{\immediate\write16{#1}}

\def\findem{\vrule height0pt width4pt depth4pt}

\long\def\demA#1{{\parindent=0pt\messages{debut de preuve}\smallbreak
     \advance\margeg by 2truecm \leftskip=\margeg  plus 0pt
     {\everypar{\leftskip =\margeg  plus 0pt}
              \everydisplay{\displaywidth=\hsize
              \advance\displaywidth  by -1truecm
              \displayindent= 1truecm}
     {\bf Proof } -- \enspace #1
      \hfill\findem}\bigbreak}\messages{fin de preuve}}

\def\resp{\mathop{\rm resp}\nolimits}
\def\resp.{\mathop{\rm resp.}\nolimits}

%FIN DES MACROS

%%%%%%%%%%%%%%%%%%%%%%%%%%%%%%%%%%%%%%%%%%%%%%%%%%%%%%%%%
%%%%%%%%%%%%%%%%%%%%%%%%%%%%%%%%%%%%%%%%%
%%%%%%%%%%%%%
\null\vskip 1cm
\centerline{\bf  ELLIPTIC CLIFFORDIAN FUNCTIONS}
\bigskip
\centerline{\bf by}
\medskip
\centerline{{\bf Guy Laville } and {\bf Ivan Ramadanoff}}

\vskip 1,5cm
{\parindent=1cm\narrower\noindent{\bf Abstract.-} \ 
 {In the study of holomorphic functions of one complex variable, one
well-known theory is that of elliptic functions and it is possible to take the
\ $\zeta$-function of Weierstrass as a building stone of this vast theory.  We
are working the analogue theory in the natural context of higher dimensional
spaces~: holomorphic and elliptic Cliffordian functions. }
\par}

\parindent=0cm
\vskip1,5cm
{\bf \SS 1. Elliptic Cliffordian functions}
\bigskip\medskip
 Everywhere in this paper we will be interested in functions \ $f
: \Omega
\fle \R_{0,3}$, \ where \ $\Omega$ \ is an open subset of \ $\R^4$.  In 
previous papers [8,9] we introduced the notion of (left) holomorphic
Cliffordian function acting from an open subset of \ $\R^{2m+2}$ \ in the
Clifford algebra \ $\R_{0,2m+1}, m\in\N$. 

\bigskip
 Let us recall the basic definition and some general properties of
holomorphic Cliffordian functions. Consider the Clifford algebra \
$\R_{0,2m+1}$ \ of the real vector space \ $V$ \ of dimension \ $2m+1$, 
provided with a quadratic form of negative signature. Denote by \ $S$ \ the
set of the scalars in \ $\R_{0,2m+1}$ \ which can be identified to \ $\R$. 
Let  \ $\{ e_i\}$, \ $i = 1,2,\ldots , 2m+1$ \ be an orthonormal basis of \
$V$ \ and let \ $e_0 = 1$. Thus, \ $\{ e_0, e_1,\ldots , e_{2m+1}\}$ \ would
be a basis for the algebra \ $\R_{0, 2m+1}$ \ such that \ $e_ie_j + e_je_i =
- 2\delta_{ij}$ \ for \ $0\leq i, j \leq 2m+1$ \ and \ $\delta_{ij}$ \ is
the Kronecker symbol.

\bigskip\bigskip
{\bf Definition.-} \ Given an open subset \ $\Omega$ \ of \ $S\oplus V$, \ a
function \ $f : \Omega\fle \R_{0,2m+1}$ \ is said to be (left) holomorphic
Cliffordian on \ $\Omega$ \ if and only if~:
$$D\Delta^m f(x) = 0$$
for each \ $x$ \ of \ $\Omega$. Here \ $\Delta^m$ \ means the iterated \ $m$ \
times Laplacian \ $\Delta$ \ and \ $D$ \ is the usual operator for
monogenic functions [1], \  [3], \  [4] :
$$D = \sum_{i=0}^{2m+1} \ e_i \ {\partial\over \partial x_i}.$$

\bigskip
 It is clear that, in the case when \ $m = 0$, such a function
would be an usual holomorphic function of one complex variable.

\bigskip
 Looking at the previous definition, one could remark that,
because of the appearance of the Laplacian, a kind of qualitative jump takes
place when the integer \ $m$ \ passes from \ $0$ \ to 1. In fact, there would
be no major differences in the studies of holomorphic Cliffordian functions
when \ $m
\geq 1$. That is the reason, and also in order to have lighter written
expressions of our formulas, we will restrict us to the case $m=1$.

 Referring again to [8,9],  let us recall that holomorphic
Cliffordian functions verify an integral representation formula which can be
looked as an analogous of the Cauchy integral formula for holomorphic
functions. Thus, we get~:

\bigskip
{\bf Theorem.}  (Taylor expansion).- \ Let \ $f$ \ be a holomorphic
Cliffordian function in an open neighborhood \ $W$ \ of \ $a\in S \oplus
V$.  Then, for each \ $x\in W$~:
$$f(x) = \sum_{k=1}^\infty \ \sum_{\vmid{\alpha} = k} \ P_\alpha (x-a)
c_\alpha ,$$

were \ $\alpha = (\alpha_0, \alpha_1, \alpha_2, \alpha_3)$ \ is a multiindex
of \ $\N^4$, \ $\vmid{\alpha}$ \ denotes its length  $\vmid{\alpha} =
\alpha_0 + \alpha_1 + \alpha_2 + \alpha_3$, \ $c_\alpha \in \R_{0,3}$ \ and
the \ $P_\alpha$ \ are left (and right) holomorphic Cliffordian polynomials
of degree \ $\vmid{\alpha} -\  1$ \ on \ $x$ \ defined as follows~: consider
the set \ $\{ e_\nu\} = \{ e_0,\ldots , e_0,\  e_1,\ldots , e_1,\ 
e_2,\ldots , e_2, \  e_3,\ldots , e_3\}$ \ where \ $e_0$ \ is written \
$\alpha_0$ \  times,
 $e_1 : \alpha_1$ \ times, \   etc $\ldots$

Set :
$$P_\alpha (x) = \sum_{\hbox{\got S}} \ \prod_{\nu = 1}^{\vmid{\alpha}-1}
(e_{\sigma (\nu )} x) \ e_{\sigma (\vmid{\alpha})}$$

the sum being expanded over all distinguishable elements \ $\sigma$ \ of the
permutation group \ \got S \ of the set \ $\{ e_\nu \}$.

\bigskip
 It should be noted that, for a fixed  \ $k\in\N^*$,  the set \
$\{ P_\alpha : \ \vmid{\alpha} = k\}$ \ contains \ $C_{k+3}^3$ \
polynomials. One could also remark that there is a generating function for
the polynomials  \
$P_\alpha$.  Indeed, if we set~:
$$\lambda = \sum_{i=0}^4 \ \lambda_i e_i, \quad \lambda_i\in\R
\quad\hbox{and}\quad \lambda_\alpha = \prod_{i=0}^4
\lambda_i^{\alpha_{i}},$$

then the following formal series gives a generating function~:
$$(1-\lambda x)^{-1} \lambda = \sum_\alpha \ P_\alpha (x) \lambda_\alpha.$$

 Using similar notations, set \ $\beta\in\N^4$ \ a multiindex of
length \
$\vmid{\beta}$.  Consider the same set \ $\{ e_\nu\}$ \ as before and put~:
$$S_\beta (x) = \sum_{\hbox{\got S}} \ \prod_{\nu = 1}^{\vmid{\beta}}
(x^{-1} e_{\sigma (\nu )}) x^{-1}.$$

Thus we will have for example \ $S_{(1,0,0,0)}(x) = x^{-1} e_0 x^{-1}$, \
$S_{(0,1,1,0)} (x) = x^{-1} e_1 x^{-1} e_2 x^{-1} + x^{-1} e_2 x^{-1} e_1
x^{-1}$.  In the rational function \ $S_\beta$ \ the power of \ $x^{-1}$ \
is \ $\vmid{\beta} + 1$.  It is not difficult to observe that the \
$S_\beta$ \ are generated by \ $x^{-1}$ \ through derivation and thus the \
$S_\beta$ \ are left (and right) holomorphic Cliffordian function excepting
at $0$. This allows us to deduce~:

\bigskip
{\bf Theorem.} (Laurent expansion).- \ If \ $B$ \ is a ball in \ $S\oplus V$
\ centered at the origin and if  \ $f : B\setminus \{ 0\} \fle \R_{0,3}$ \ is
a holomorphic Cliffordian function, then for each \ $x\in B \setminus \{
0\}$~:
$$f(x) = \sum_{\vmid{\beta} = 0}^\infty S_\beta (x) d_\beta +
\sum_{\vmid{\alpha} = 1}^\infty P_\alpha (x) c_\alpha ,$$
where \ $c_\alpha$ \ and \ $d_\beta$ \ belong to \ $\R_{0,3}$.

\bigskip
 In such a way, the first sum should be considered as the
analogous of the singular part of a Laurent expansion for a holomorphic
function, while the second sum represents the analogous of its regular part.

\bigskip
Let us make some additional remarks.  If \ $f$ \ is a real
analytic function on a neighborhood \ $W$ \ of \ $a\in S \oplus W$ \ and
taking its values in \ $S\oplus V$,  we can look at \ $f$ \ as a \
$C^\infty$ \ application of \ $W\subset \R^4$ \ into \ $\R^4$ \ so that \
$f$ \ possesses an usual Taylor expansion which can be written as~:
$$f(a+h) = \sum_{n=0}^\infty \ {1\over n!} \ \bigl( h \mid \nabla_x)^n
f(x)\Big\vert{_{x=a}},$$

where \  \ $( h \mid \nabla_x)$ \ is  the usual scalar product in \
$\R^4$ \ of the two vectors \ $h = (h_0, h_1, h_2, h_3)$ \ and the gradient
\ $\nabla_x = \Bigl( \displaystyle{\partial\over \partial x_0},
\displaystyle{\partial\over\partial x_1},
\displaystyle{\partial\over\partial x_2},
\displaystyle{\partial\over\partial x_3}\Bigr)$ \ applied on \ $f$ \ as a
function of  $x = (x_0, x_1, x_2, x_3)$~; \ $(h \mid \nabla )^n$
means \ $(h \mid \nabla ) \cdots (h \mid \nabla )$ \ repeated \
$n$ \ times and, of course \ $a+h \in S \oplus V$.

They are some formal calculations related to the operator \ $(h \mid
\nabla_x)$ \ which should be noted. For instance~: $(h \mid
\nabla_x) (x) = h$ \ and since \ $xx^{-1} = 1$, then  \ $(h \mid
\nabla_x) (xx^{-1}) = 0$,  but also
$$h x^{-1} + x (h \mid \nabla_x) (x^{-1}) = 0$$
so that :
$$(h \mid \nabla_x) (x^{-1}) = - x^{-1} h x^{-1}.$$

By iteration on \ $q\in\N$, \ one get :
$$(h \mid \nabla_x)^q (x^{-1}) = (-1)^q q~! \ (x^{-1} h)^q x^{-1}.$$

\bigskip\medskip
 Remember that the statement of the classical theorem of
Liouville for holomorphic functions remains valid for harmonic functions and
also for monogenic functions ([1], p.99). This allows us to state the
following~:

\bigskip\medskip
\noindent{\bf Theorem of Liouville.-} \  If \ $f$ \ is a
holomorphic Cliffordian function in the whole open set \ $S\oplus V$,
identified to \ $\R^4$, \ (such a function could be called
\undertext{entire} Cliffordian)  and if \ $f$ \ and its derivatives up to the
second order are bounded, then \ $f$ \ reduces to a constant.

\bigskip\medskip
 Whitout loss of generality we may achieve the proof for all the
components of \ $f$ \ taking its values in \ $\R_{0,3}$, \ so in what
follows we may assume that \ $f : \R^4 \fle \R$.  Set \ $g = \Delta f$. Thus
\ $g$ \ is monogenic and the Liouville theorem for monogenic functions says
that \ $g$ \ reduces to a constant, \ie \ $\Delta f = K$, \ $K\in\R$. 
Remark the condition \ $g$ \ bounded is fulffiled because of the boundedness
of the derivatives of $f$.   But \ $\Delta f = K$ \ implies \ $\Delta \bigl(
\displaystyle{\partial f\over \partial x_i}\bigr) = 0$, \ $i = 0,1,2,3$.
Apply now the Liouville theorem for harmonic functions, one deduces \
$\displaystyle{\partial f\over \partial x_i} = c_i$, \ $i = 0,1,2,3$, \
$c_i\in\R$.  This means that
$$f(x_0, x_1, x_2, x_3) - \sum_{i=0}^3 c_ix_i$$
satisfies \ $df = 0$, \ from which we get that \ $f$ \ reduces to a
polynomial on \ $x_0,x_1,x_2,x_3$ \ of first degree. Remember again \ $f$ \
must be bounded, so that \ $f$ \ reduces to a constant.

\bigskip\bigskip
 It should be noted that the same argument allows to get a
Liouville theorem for biharmonic functions under the assumption such a
biharmonic function should be bounded with all its derivatives up to the
second order.

\bigskip
 Remark that all we said on the holomorphic Cliffordian
functions, including their Laurent expansions, allows us to consider
meromorphic Cliffordian functions in the sense that such a function would be
holomorphic Cliffordian in an open subset \ $\Omega$ \ of \ $S\oplus V$,
excepting at isolated pointwise singularities, namely pointwise poles, \ie
such that the principal part of their Laurent expansion contains only a
finite number of terms. So, put~:

\bigskip\bigskip
{\bf Definition.-} \ A function \ $f : \Omega\fle\R_{0,3}$, \ where \
$\Omega$ \ is an open subset of \ $S\oplus V$ \ is said to be meromorphic
Cliffordian if \ $f$ \ is holomorphic Cliffordian on $\Omega$, excepting on
an isolated set of pointwise poles.

\bigskip
 In this sense \ $x^{-1}$ \ is a meromorphic function in \
$S\oplus V$ \ with a simple pole at the origin.

\bigskip
 The variety of singularities of general meromorphic functions is
larger. As the following example shows, the set of singularities of a
Cliffordian function might be quite complicated~:
\medskip
 Look at \ $g : x\longmapsto (x_1 - e_1x_0)^{-1}$, \ which is
holomorphic Cliffordian everywhere except at the subset \ $\{ x_1 = 0, \ x_0
= 0\}$ \ of \ $\R^4$.

\bigskip
 The study of the singularities sets which are not pointwise
deserves to  be accomplished but will not be the subject of this paper.

\bigskip
 We will end this paragraph putting the basis of elliptic
Cliffordian function theory. As in the case of functions of a complex
variable, where they are two kinds of periodic functions -- the 1-periodic
and the doubly periodic functions, here we could expect much more subclasses
of periodic functions, namely the $N$-periodic functions, where \ $N =
1,2,3,4$.

\medskip
 So take \ $N$ \ in \ $\{ 1,2,3,4\}$.  Let \ $\omega_\alpha\in
S\oplus V$ \ for \ $\alpha = 1,2,\ldots ,N$.  The \ $\omega_\alpha$ \ will
play the role of half periods. We will always suppose the paravectors \
$\omega_1,\ldots,\omega_N$ \ linearly independant in \ $S\oplus V$.

\bigskip\bigskip
{\bf Definition.-} \ A function \ $f : S\oplus V \fle \R_{0,3}$ \ is said to
be \ $N$-periodic if
$$f \bigl( x + 2  \omega_\alpha\bigr) = f(x),$$
\vskip-0,5cm
for every \ $x\in S\oplus V$.

\bigskip\medskip
 Introduce \ $\omega = (\omega_1,\ldots,\omega_N)$ \ and, for a
multiindex \
$k = (k_1,\ldots,k_N)$ \ belonging to \ $\Z^N$, \ define \ $k\omega$ \ by \
$k\omega = \displaystyle\sum_{\alpha = 1}^N k_\alpha \omega_\alpha$.  Let us
call the set  \ $2\Z^N\omega = \{ 2k\omega , k\in\Z^N\}$ \ a lattice. A
fundamental cell of this lattice will be the hypercube spanned over \
$2\omega_1$, $2\omega_2$, $2\omega_3$, $2\omega_4$.

\bigskip
 Obviously, for a \ $N$-periodic function $f$, we have :
$$f(x+ 2k\omega ) = f(x).$$
Now, define :
\bigskip
{\bf Definition.-} \ A function \ $f : S\oplus V \fle \R_{0,3}$ \ is said
elliptic Cliffordian if \ $f$ \ is \ $N$-periodic and meromorphic.

\bigskip\medskip
 The first fundamental result of the classical theory of
elliptic functions says that an elliptic function which is holomorphic
reduces to a constant. Here we have the same~:

\bigskip\medskip
{\bf Theorem 1.-} \ An elliptic Cliffordian function which is holomorphic
Cliffordian reduces to a constant.

\bigskip\bigskip
 Indeed, a holomorphic Cliffordian function is obviously bounded
on a fundamental cell and by periodicity would be bounded on the whole
espace \ $S\oplus V$.  The same is true for all its derivatives because if \
$f$ \ is \ $N$-periodic,  then all the \ $(e_i \mid \nabla_x) f(x)$ \
are also $N$-periodic. All of them would be bounded on
\
$S\oplus V$
\ via the fundamental cell. It remains just to apply the Liouville's theorem
and to conclude that \ $f$ \ reduces to a constant.

\bigskip\bigskip
 It should be remarked that the present theory of elliptic
Cliffordian function is only an additive theory in the sense that elliptic
Cliffordian functions could be added and multiplied (in \ $\R_{0,3}$)  by
Clifford constants.

\bigskip
 At this time we do not dispose with a satisfactory definition
of a product of two elliptic Cliffordian functions.

\bigskip
 The next theorem is an example of what we can do only with
additive arguments~:

\bigskip
{\bf Theorem 2.-} \ If \ $f_1$ \ and \ $f_2$ \ are two elliptic Cliffordian
functions with the same pointwise poles and the same principal part of
their Laurent expansions on the neighborhoods of their poles, then they
differ just up to an additive constant.

\bigskip\medskip
 The proof is easy : consider \ $f_1-f_2$.  This is a
holomorphic Cliffordian function which is periodic and by theorem 1 must be a
constant.

\vskip 1cm
{\bf \SSg\kern.15em  2. On the trigonometric Cliffordian functions}
\bigskip\medskip
\font\grandsy=cmsy10 scaled\magstep0
\def\SSg{{\grandsy x}}
 Let us start with a slight modification of the lemma of
\SSg 4 in [9] \ which says~:

\bigskip\medskip
\noindent{\pmb{\sc Lemma 1.-}} \ {\it If \ $u : \R^2 \fle \R , \ (\xi , \eta
) \longmapsto u (\xi , \eta )$ \ is harmonic, then \ ${\cal U} (x) = u(x_0,
\vmid{\hskip-0,10cm\vec x \hskip-0,10cm})$, \ where \ $x = x_0 + \vec x$, \
is biharmonic, \ i.e. \ $\Delta^2 {\cal U} (x) = 0$. 
}

\bigskip
 With the help of this lemma, we dispose a way to generate
holomorphic Cliffordian functions, taking \ $D^* {\cal U}$. Recall that~:
$$D = {\partial\over\partial x_0} + \sum_{i=1}^3 \ e_i \
{\partial\over\partial x_i}$$
and \ $D^* = \displaystyle{\partial\over\partial x_0} -
\displaystyle\sum_{i=1}^3 \ e_i \ \displaystyle{\partial\over\partial x_i}$,
\ so that if we introduce the notations \ $D_0 = \displaystyle{\partial\over
\partial x_0}$ \ and \ $\overrightarrow D = \displaystyle\sum_{i=1}^3 \ e_i
\ \displaystyle{\partial\over \partial x_i}$, \ we have \ $D = D_0 +
\overrightarrow D$ \ and \ $D^* = D_0 - \overrightarrow D$.

\bigskip\bigskip\medskip
 Now, let us modify the above lemma in order to have a way for
generating holomorphic Cliffordian functions directly from the holomorphic
ones~:

\bigskip\bigskip
\noindent{\pmb{\sc Lemma 2.-}} \ {\it  If \ $f : (\xi , \eta ) \longmapsto f
= u+iv$ \ is a holomorphic function, then  $F(x) = u(x_0, \ \vmid{\vec
x}) + \displaystyle{\vec x\over \vmid{\vec x}} v (x_0, \ \vmid{\vec x})$ \ is
a holomorphic Cliffordian function.
} 
\bigskip
{\bf Proof.-} Given \ $f$ \ holomorphic, there exists $h$, which is
harmonic and \ $2 \ \displaystyle{\partial\over \partial z} h = f$, \ so that 
$\displaystyle{\partial h\over \partial \xi} = u$, \ $\displaystyle{\partial
h\over \partial \eta} = - v$. Apply lemma 1 to $h$, we construct \ $H(x) = h
(x_0, \ \vmid{\vec x})$ \ which will be biharmonic and thus \ $F(x) = D^*
H(x)$ \ should be holomorphic Cliffordian. But~:
$$D^* H(x) = D_0 H(x) - \overrightarrow D H(x).$$

However,
$$\eqalign{
&D_0H (x) = {\partial h\over \partial x_0} (x_0, \vmid{\vec x}) = {\partial
h\over \partial \xi} (\xi , \eta ) \Bigl|_{\scriptstyle \xi = x_{0},\eta =
\vmid{\vec x}} = u(x_0, \vmid{\vec x}), \ \hbox{ and }\cr
&\overrightarrow D H(x) = {\partial h\over \partial\eta} (\xi , \eta
)\Big|_{\xi = x_{0}, \eta = \ \vmid{\vec x}} \cdot \overrightarrow D
(\vmid{\vec x}), \quad \hbox{and}\cr
&\overrightarrow D (\vmid{\vec x}) = \overrightarrow D {\bigl( (\vmid{\vec
x}^2\bigr)}^{1\over 2} = {1\over 2} \ {1\over \vmid{\vec x}} \
\overrightarrow D (\vmid{\vec x}^2) = {\vec x\over \vmid{\vec x}},\cr}$$ 

so that \ $\overrightarrow D H(x) = - \displaystyle{\vec x\over \vmid{\vec
x}} v (x_0, \ \vmid{\vec x})$ \ and finally~:
$$D^* H(x) = F(x) = u (x_0, \ \vmid{\vec x}) + {\vec x\over \vmid{\vec x}} v
(x_0, \ \vmid{\vec x}).$$

\bigskip\bigskip
 Lemma 2 is a fruitful tool for generating  the analogous of
the trigonometric functions \ $\sin z$, $\cos z$ etc $\ldots$ \ in the 
frame of Cliffordian analysis.

\bigskip
 Thus, for example, the function
$$e^x = e^{x_0} \ \bigl( \cos \sqrt{-\vec x^2} + {\vec x\over \vmid{\vec x}}
\ \sin \sqrt{- \vec x^2}\bigr)$$

is a well-defined holomorphic Cliffordian function from \ $S\oplus V$ \ in \
$S\oplus V \subset \R_{0,3}$.

\bigskip
 Further,
$$\eqalign{
&\sin x = - {\vec x\over \vmid{\vec x}} \ {e^x - e^{-x}\over 2} , \cr
&\cos x = {e^x + e^{-x}\over 2} , \cr
&\cotan x = {\vec x\over \vmid{\vec x}} \ {e^x + e^{-x}\over e^x - e^{-x}}
.\cr}$$

\bigskip
 Lemma 2 allows us to consider all of them as holomorphic
Cliffordian functions. The final procedure for constructing such functions
is to make a formal change of the letters \ $z$ \ in $x$, \ $\xi = {\cal R}e
\{ z\}$ \ in \ $x_0$, \ $i\eta = i \Im \{ z\}$ \ in \ $\vec x$, or
equivalently \ $i$ \ in \ $\displaystyle{\vec x\over \vmid{\vec x}}$ \ and \
$\eta$ \ in \ $\vmid{\vec x}$.

\vskip1cm
{\bf 3. Weierstrass \ $\zeta$ \ holomorphic Cliffordian functions}
\bigskip\medskip
  Remember that the point of view of Weierstrass theory of elliptic
functions (doubly-periodic and meromorphic) is based on the ${\cal P}$-function
which appears naturally as the derivative, up to the sign, of the \ $\zeta$-
function of Weierstrass.

\bigskip
 In this paragraph we will construct  four \ $\zeta_N$ \
Weierstrass functions in \ $\R_{0,3}$. Take again
\ $N$ \ in
\ $\{ 1,2,3,4\}$.  Introduce \ $\omega = (\omega_1, \ldots ,\omega_N)$ \
where \ $\omega_\alpha \in S \oplus V$ \ for  $\alpha = 1,2,\ldots
,N$,   and suppose the paravectors \ $\omega_1,\ldots , \omega_N$ \ linearly
independant in \ $S\oplus V$. For a multiindex \ $k = (k_1,\ldots , k_N)$ \
belonging to \ $\Z^N$,  define  \ $k\omega$ \ by \ $k\omega =
\displaystyle\sum_{\alpha = 1}^N k_\alpha \omega_\alpha$. Taking into
account the countability of the respective sets, it is not difficult to
rearrange the lattice \ $2 \Z^N \omega \setminus \{ (0,\ldots,0)\}$ \ as \
${\{ w_p\}}_{p=1}^\infty$.

\bigskip\bigskip
{\bf Definition.-} \ A Weierstrass \ $\zeta_N$-function with lattice \
$2\Z^N \omega$ \ is the function  $\zeta_N : S \oplus V \setminus 2\Z^N
\omega \fle \R_{0,3}$ \  defined by :
$$\zeta_N(x) = x^{-1} + \sum_{p=1}^\infty \ \left\{ (x- w_p)^{-1} + \sum_{\mu
= 0}^{N-1} (w_p^{-1}x)^\mu \ w_p^{-1}\right\} .\leqno (1)$$

\hskip 0,45cm In such a way we get four functions \ $\zeta_1, \zeta_2,
\zeta_3, \zeta_4$
\ which are related respectively with the 1-periodical lattice \ $\{ 2k_1
\omega_1 , k_1\in \Z , \omega_1\in S \oplus V\}$ \ for the first one and
with a 4-periodical lattice for the last one.

\bigskip
\hskip 0,45cm Consider a particular case concerning \ $\zeta_1$, \ setting \
$\omega_1 =
\displaystyle{\pi\over 2}$~:
$$\zeta_1(x) = x^{-1} + \sum_{k_{1}\in \Z^*} \ \{ (x-k_1\pi )^{-1} + (k_1\pi
)^{-1}\} .$$

Thus we recognize the analogous of the function \ $\cotg z$ \ in the function
theory of one complex variable.

\bigskip
\hskip 0,45cm Remark also that \ $\zeta_2$ \ is analogous to the well-known
Weierstrass \
$\zeta$-function of a complex variable~:
$$\zeta (z) = {1\over z} + \sum_{p=1}^\infty \ \left\{ {1\over z-w_p} +
{1\over w_p} + {z\over w_p^2} \right\} .$$

Here obviously, the subset \ $(2\Z\omega_1 + 2\Z\omega_2) \setminus \{
(0,0)\}$ \ of the corresponding lattice has been rearranged in a set noted \
$\{ w_p\},\  p\in\N^*$.

\bigskip\bigskip
 So we have obtained four \ $\zeta$-functions. If we pursue in the
way of the similarity with holomorphic functions, we could expect that all of
them be holomorphic Cliffordian excepted at the points of the lattice. They
probably would not be themselves \ $N$-periodic, but by \ $N$ \ derivations
they would generated the analogous \ ${\cal P}_N$ \ of the Weierstrass \
${\cal P}$-function.

First of all, note that the series
$$\sum_{w\in\Omega_{N}, w\not= 0} \ \vmid{w}^{-N-1}$$
converges where \ $\Omega_N$ \ designes a \ $N$-periodic lattice in \
$\R^4$.  This is a result following easily from similar arguments as in
the classical case.

\bigskip
 The existences of the \ $\zeta_N$ \ will be deduced from the
proof of the convergences of the series defining \ $\zeta_N$.  For this, a
majoration of the general term of the series will be achieved. First, remark
the general term of the series could be written as~:
$$\displaylines{
(x-w_p)^{-1} + \sum_{\mu = 0}^{N-1} (w_p^{-1} x)^\mu w_p^{-1} = \left\{ 1 +
\sum_{\mu = 0}^{N-1} (w_p^{-1}x)^\mu w_p^{-1} (x-w_p)\right\} (x-w_p)^{-1}
=\cr
= \left\{ \sum_{\nu = 0}^N {(w_p^{-1} x)}^\nu - \sum_{\mu = 0}^{N-1}
{(w_p^{-1} x)}^\mu\right\} (x-w_p)^{-1} = (w_p^{-1} x)^N (x-w_p)^{-1}.\cr}$$

\medskip
Now, let us consider the compact set \ $K_r = \{ x \in S \oplus V : \
\vmid{x}
\ \leq r\}$.  Thus, except for a finite number of \ $w_p$, one has \
$\vmid{w_p} \ \geq 2r$. 

For \ $\vmid{x} \ \leq r$, \ the next majoration follows :
$$\displaylines{\vmid{(w_p^{-1}x)^N (x-w_p)^{-1}} \  = \ \vmid{
{(w_p^{-1}x)}^N {(w_p^{-1} x-1)}^{-1} w_p^{-1}} \  \leq\cr
\leq {r^N\over 2} \ {\vmid{w_p^{-1}}}^{N+1} = C_{N,r}
\vmid{w_p}^{-N-1}.\cr}$$

In such a way we proved that the series defining \ $\zeta_N$ \ are uniformly
convergent on compact subsets of \ $S\oplus V$ \ for \ $N = 1,2,3,4$.

\bigskip
 It should be noted that the limit (or the sum) of every
sequence (or series) of holomorphic Cliffordian functions converging
uniformly on compact sets is also a holomorphic Cliffordian function. The
reason is that \ $f$ \ is holomorphic Cliffordian if and only if \ $DD^*Df =
0$,  so that applying three times the corresponding theorem of Weierstrass
for monogenic functions [1], p.58, we get the result.

\bigskip
 In our case, all the terms in the series defining \ $\zeta_N$ \
are holomorphic Cliffordian in \ $S\oplus V$ \ excepting at the points \
$w_p$ \ and the convergence being uniform on compact subsets, we can say
that \ $\zeta_N$ \ is a holomorphic Cliffordian function on \ $S\oplus V
\setminus {\{ w_p\}}_{p=0}^\infty$, \ where \ $w_0 = (0,0,\ldots, 0)$, \ or
equivalently on \ $S\oplus V \setminus 2\Z^N \omega$.

\bigskip
 Generally, a holomorphic Cliffordian function acts from \
$S\oplus V$ \  in \ $\R_{0,3}$.  As far as the \  $\zeta_N$ \ functions are
concerned, the first observation we can do, is that \ $\zeta_N$ \ possesses
simple poles at the points of the lattice (in fact, by construction).  We
could remark also that \ $\zeta_N$ \ takes its values not in the whole
Clifford algebra \ $\R_{0,3}$, \ but in a subset, which is exactly \
$S\oplus V$.

\bigskip
 Indeed, \ $x^{-1} = \displaystyle{x^*\over \vmid{x}^2}$, \
where \ $x^*$ \ is the conjugate of $x$, \ $x^* = x_0e_0 -
\displaystyle\sum_{j=1}^3 x_j e_j$, obviously belongs to \ $S\oplus V$. A
straightforward computation carried on \ $hxh$, \ where \ $x,h\in S\oplus V$
\ shows that \  $hxh \in S \oplus V$ \ and then, using a reccurrence
argument, one can show that \ $(hx)^n h \in S\oplus V$ \ for \ $n\in\N$.

 At the end of this paragraph let us see that the \ $\zeta_N$ \
functions are odd. In fact, there is another way to write down the
expressions of the \ $\zeta_N$ \ functions, namely~:
$$\zeta_N(x) = x^{-1} - \sum_{n=N}^\infty \ \sum_{p=1}^\infty  (w_p^{-1}
x)^n w_p^{-1} \leqno (2)$$

which is easily obtained developping \ $(x-w_p)^{-1}$ \ as \ $-w_p^{-1} -
w_p^{-1} x w_p^{-1} - \cdots -  (w_p^{-1} x)^n w_p^{-1} - \cdots $ \ in
(1).  Here, we have to remember that the set \ ${\{ w_p\}}_{p=1}^\infty$ \
is exactly \ $2\Z^N\omega \setminus \{ (0,\ldots ,0)\}$.  Therefore, the sum
(2) contains terms of the type
$${\{ (2k\omega )^{-1}x (2k \omega )^{-1}\}}^n (2k\omega )^{-1}$$

for \ $k\in\Z^N \setminus \{ (0,\ldots , 0)\}$ \ and \ $n\geq N$.  But when
\ $n$ \ is even, we meet  ``opposite" vertices of the lattice, i.e. such
that the lengths of the corresponding multiindices are  $+ \vmid{k}$ \ and \
$- \vmid{k}$ \ and thus the corresponding terms in the sum will cancelled.
So the sum is expanded over those \ $n$ \ which are odd integers and we get 
$$\zeta_N (-x) = - \zeta_N (x)$$
for \ $x\in S \oplus V \setminus {\{ w_p\}}_{p=0}^\infty$ \ and \ $N =
1,2,3,4$.  The formula (2) is reduced to~:
$$\zeta_N(x) = x^{-1} - \sum_{k \geq [{N\over 2}]} \ \sum_{p=1}^\infty
(w_p^{-1}x)^{2k+1} w_p^{-1}. \leqno (3)$$

Remark this formula co\"\i ncides with the well know formula for the classical
Weierstrass \ $\zeta$ \ function in the case $N=2$.  That we obtained is~:
$$\displaylines{ \zeta_2(x) = x^{-1} - \sum_{k\geq 1} \ \sum_{p=1}^\infty \
(w_p^{-1}x)^{2k+1} w_p^{-1} =\cr
= x^{-1} - \sum_{p=1}^\infty (w_p^{-1} x)^3 w_p^{-1} - \sum_{p=1}^\infty
(w_p^{-1} x)^5 w_p^{-1} - \cdots \cr}$$

which is the analogous of :
$$\zeta (z) = {1\over z} - z^3 \ \sum_{p=1}^\infty {1\over w_p^4} - z^5 \
\sum_{p=1}^\infty {1\over w_p^6} - \cdots \leqno (3')$$

\bigskip
 We will end this paragraph by the Laurent expansion of \
$\zeta_N$ \ in a neighborhood of the arign. This can be achieved shortly
combining (2) with the last formula of \ \SS 1 we could write as~:
$$(x \mid \nabla )^n (w_p^{-1}) = (-1)^n n! \ (w_p^{-1} x)^n \
w_p^{-1},$$

where \ $x, w_p\in S \oplus V$~; \ $p,n\in\N$.  So, we get~:
$$\eqalign{ \zeta_N(x)
&= x^{-1} - \sum_{n=N}^\infty \ \sum_{p=1}^\infty (-1)^n \ {(x \mid
\nabla )^n\over n!} (w_p^{-1}) \cr
&= x^{-1} - \sum_{n=N}^\infty (-1)^n \ {(x \mid \nabla )^n\over n!} \
\Bigl( \sum_{p=1}^\infty w_p^{-1}\Bigr) .\cr}$$

Taking into account (3), we have :
$$\zeta_N(x) = x^{-1} + \sum_{k \geq [{N\over 2}]} \ {(x \mid \nabla
)^{2k+1}\over (2k+1)!} \ \Bigl( \sum_{p=1}^\infty w_p^{-1}\Bigr) . \leqno
(4)$$

In this way we are able to give another form to the classical formula (3'),
namely~:
$$\zeta (z) = {1\over z} + {(z \mid \nabla )^3\over 3!} \ \Bigl(
\sum_{p=1}^\infty {1\over w_p}\Bigr) +  {(z \mid \nabla )^5\over 5!}
\  \Bigl( \sum_{p=1}^\infty {1\over w_p}\Bigr) + \cdots$$

\vskip1cm
{\bf 4. Quasi-periodicity of the \ $\zeta_N$ \ functions}
\bigskip\medskip
 The aim of this paragraph is to obtain the analogous of the
formula of quasi-periodicity of the classical Weierstrass function in the
case of a complex variable. It is well-known that there exist two constants
\ $\eta , \eta'$, \ such that~:
$$\zeta (z+2n\omega + 2m\omega') = \zeta (z) + 2n\eta + 2m\eta'$$

where \ $z\in \C , n, m\in\Z$,  \ and the corresponding lattice is generated
by \ $2\omega$ \ and \ $2\omega'$.  Moreover, \ $\eta = \zeta (\omega )$ \ and
\ $\eta' = \zeta (\omega')$.
The first step is to show that the function

$Z_4(x,a) = \zeta_4 (x+a) - \zeta_4(a) - (x\mid \nabla ) \zeta_4 (a)
- \displaystyle{1\over 2!} (x \mid \nabla )^2 \zeta_4 (a) -
\displaystyle{1\over 3!} (x \mid \nabla )^3 \zeta_4 (a)$  \  looked
as a function of the variable \ $a$ \  is holomorphic Cliffordian on \
$S\oplus V
\setminus {\{ w_m\}}_{m=0}^\infty$ \ and also is periodic. In fact, we have
to substitute all the terms of this sum by their expansions and observe that
all the terms containing polynomial expressions on \ $a$ \  are cancelled. The
remaining part will satisfy the two required conditions. The computations
will be carried on \ $\zeta_4$ \ denoted for briefness by  $\zeta$.

\medskip
 A straightforward calculus gives~:
$$\eqalign{\zeta (x+a) 
&= (x+a)^{-1} + \sum_{m=1}^\infty \{ (x+a-w_m)^{-1} + w_m^{-1} + w_m^{-1}
(x+a) w_m^{-1}\cr
&+ (w_m^{-1} \bigl( x+a)\bigr)^2 w_m^{-1} + \bigl( w_m^{-1} \bigl(
x+a)\bigr)^3 w_m^{-1} \} \cr}$$

which can be written as : 
$$\eqalign{
&\zeta (x+a) = \sum_{m=0}^\infty (x+a-w_m)^{-1} + \sum_{m=1}^\infty \ \Bigl\{
w_m^{-1} + w_m^{-1} x \ w_m^{-1} + w_m^{-1} a \ w_m^{-1} +\cr
&+ {(w_m^{-1} x)}^2 w_m^{-1} + w_m^{-1}x\  w_m^{-1} a\  w_m^{-1} + w_m^{-1}
a\  w_m^{-1} x \ w_m^{-1} + {(w_m^{-1}a)}^2 w_m^{-1}+\cr
&+ {(w_m^{-1} x)}^3 w_m^{-1} + {(w_m^{-1} x)}^2 w_m^{-1} a \ w_m^{-1} +
w_m^{-1} x \ w_m^{-1} a \  w_m^{-1} x \ w_m^{-1} +\cr
&+ w_m^{-1} a \ {(w_m^{-1} x)}^2 \ w_m^{-1} + {(x_m^{-1}a)}^2 w_m^{-1} x
w_m^{-1} + w_m^{-1} a \  w_m^{-1} x \  w_m^{-1} a \  w_m^{-1} +\cr
&+ w_m^{-1} x \ {(x_m^{-1} a)}^2 w_m^{-1} + {(w_m^{-1} a)}^3 w_m^{-1}  \Bigr\}
\cr}$$

where we have reintroduced in the first sum the power \ $(x+a)^{-1}$ \
appearing for \ $w_0=0$.

From this sum we have to substract four other sums :
$$\zeta (a) = \sum_{m=0}^\infty (a - w_m)^{-1} + \sum_{m=1}^\infty \ \Bigl\{
w_m^{-1} + w_m^{-1} a \  w_m^{-1} + {(w_m^{-1}a)}^2 w_m^{-1} + {(w_m^{-1}
a)}^3 w_m^{-1}\Bigr\}$$

and :
$$\eqalign{
&(x \mid \nabla ) \zeta (a) = - \sum_{m=0}^\infty (a-w_m)^{-1} x
(a-w_m)^{-1} + \sum_{m=1}^\infty \Bigl\{ w_m^{-1} x \  w_m^{-1} +\cr
&+ w_m^{-1} x \  w_m^{-1} a \  w_m^{-1} + w_m^{-1} a \  w_m^{-1} x \
w_m^{-1} + {(w_m^{-1} a)}^2 w_m^{-1} x \  w_m^{-1}\cr
&+ w_m^{-1} a \  w_m^{-1} x w_m^{-1} a \ w_m^{-1} + w_m^{-1} x \ {(w_m^{-1}
a)}^2 w_m^{-1}\Bigr\} . \cr}$$

Further :
$$\eqalign{
&{1\over 2} (x \mid \nabla )^2 \zeta (a) = \sum_{m=0}^\infty {\bigl(
(a-w_m)^{-1} x\bigr)}^2 (a-w_m)^{-1} + \sum_{m=1}^\infty \ \Bigl\{
{(w_m^{-1} x)}^2 w_m^{-1} \cr
&\qquad + {(w_m^{-1}x)}^2 w_m^{-1} a \  w_m^{-1} + w_m^{-1} x \  w_m^{-1} a
\  w_m^{-1} x \  w_m^{-1} + w_m^{-1} a \ {(w_m^{-1}x)}^2 w_m^{-1} \Bigr\}
,\cr}$$

and finally : 
$${1\over 3!} (x \mid \nabla )^3 \zeta (a) = - \sum_{m=0}^\infty \
{\bigl( (a-w_m)^{-1}x\bigr)}^3 (a-w_m)^{-1}  + \sum_{m=1}^\infty {(w_m^{-1}
x)}^3 w_m^{-1}.$$ 

A careful substraction gives :
$$\eqalign{&Z_4(x,a) = \zeta (x+a) - \zeta (a) - (x\mid \nabla ) \zeta
(a) - {1\over 2!} (x \mid \nabla )^2 \zeta (a) - {1\over 3!} (x \mid
\nabla )^3 \zeta (a) =\cr
&= \sum_{m=0}^\infty (x+a-w_m)^{-1} - \sum_{m=0}^\infty (a-w_m)^{-1}  +
\sum_{m=0}^\infty (a-w_m)^{-1} x (a-w_m)^{-1} \cr
&- \sum_{m=0}^\infty {\bigl( (a-w_m)^{-1}x\bigr)}^2 (a-w_m)^{-1} +
\sum_{m=0}^\infty {\bigl( (a-w_m)^{-1}x\bigr)}^3 (a-w_m)^{-1}.\cr}$$

\bigskip
\hskip 0,45cm If we consider \ $Z_4(x,a)$ \ as a function of the variable
$a$, because every term in each sum on the right hand side is holomorphic
Cliffordian and the convergence is fulffiled then \ $Z_4(x,a)$ \ is
holomorphic Cliffordian on \ $a$\  over \ $S\oplus V \setminus {\{
w_n\}}_{m=0}^\infty$. Remember the lattice \ ${\{ w_n\}}_{m=0}^\infty$ \ can
be seing as  \ $2\Z^4\omega$, \ all the sums are periodic and so \
$Z_4(x,a)$ \ is periodic on \  $a$.

Let us write down this periodicity :
$$Z_4 (x, a+2\omega ) = Z_4 (x,a)$$
and put here \ $a = -\omega$, \ i.e. :
$$\eqalign{
&\zeta (x+\omega ) - \zeta (\omega ) - (x\mid\nabla ) \zeta (\omega )
- {1\over 2!} (x\mid\nabla )^2 \zeta (\omega )\cr
&- {1\over 3!} (x\mid \nabla )^3 \zeta (\omega ) = \zeta (x-\omega ) -
\zeta (-\omega ) - (x\mid\nabla ) \zeta (-\omega )\cr
&- {1\over 2!} (x\mid\nabla )^2 \zeta (-\omega ) - {1\over 3!}
(x\mid\nabla )^3 \zeta (-\omega ).\cr}$$

Recall \ $\zeta$ \ is odd, so \ $(x \mid\nabla )^2\zeta$ \ will be
also odd and thus~:
$$\zeta (x+\omega ) - \zeta (x-\omega ) = 2\zeta (\omega ) + (x\mid
\nabla )^2 \zeta (\omega ).$$
Apply now the last formula for \ $x+\omega$, \ we get :
$$\zeta_4(x+2\omega ) - \zeta_4(x) = 2\zeta_4 (\omega ) + \bigl( (x+\omega
)\mid\nabla \bigr)^2 \zeta_4(\omega )\leqno (5)$$

which shows the quasi-periodicity of the \ $\zeta_4$ \ Weierstrass function.
The right hand side is a polynomial on \ $x$ \ of degree 2.

\bigskip
\hskip 0,45cm It should be noted that following the same way, we are in a
position to get quasi-periodicity formulas for all the \ $\zeta_N$ \
Weierstrass functions. For \ $\zeta_1$ \ a direct computation shows that
$$\zeta_1 (x+2\omega ) - \zeta_1(x) = 0.$$

For the other cases, \ $N = 2,3,4$, \ it suffices to erase all the terms
which should disappear in the computation of the corresponding \ $Z_N(x,a)$.
For instance~:
$$Z_N(x,a) = \zeta_N(x+a) - \sum_{n=0}^{N-1} \ {(x\mid \nabla
)^n\over n!} \zeta_N(a), \ \ N = 2,3,4.$$

Those are \ $N$-periodic. Writing down the periodicity and taking into
account that \ $\zeta (\omega )$, \ $\displaystyle{(x\mid\nabla
)^2\over 2!} \zeta (\omega )$ \ are odd on $\omega$, one get~:
$$\zeta_2 (x+2\omega ) - \zeta_2(x) = 2\zeta_2(\omega )$$
and
$$\zeta_3 (x+2\omega ) - \zeta_3(x) = 2\zeta_3(\omega ) + \bigl( (x+\omega
)\mid\nabla\bigr)^2 \zeta_3(\omega ).$$

\hskip 0,45cm Let us introduce a notation fo the right-hand side polynomial
of quasi-periodicity as, for example,
$$\zeta_N (x+2\omega ) - \zeta_N(x) = \eta_N(x,\omega ).$$

\hskip 0,45cm Except in the case \ $N = 1$, \ when \ $\eta_1$ \  reduces to
the constant zero, \ $\eta_N$ \ contains terms of the type \ $2 \
\displaystyle{\bigl( (x+\omega )\mid\nabla \bigr)^{2p}\over (2p)!} \
\zeta_N(\omega )$ \ with \ $p\in \left[\hskip-0,1cm\left[ 0,
\left[\displaystyle{N+1\over 2}\right] - 1 \right]\hskip-0,1cm\right]$ \ when
\
$N\geq 2$.  So that~:
$$\eta_N(x,\omega ) = 2 \ \sum_{p=0}^{[{N+1\over 2}] -1} \ {\bigl( (x+\omega
)\mid\nabla \bigr)^{2p}\over (2p)!} \zeta_N(\omega ).$$

\vskip 1cm
{\bf 5. Some other formulas for the \ $\zeta_4$ \ function}
\bigskip\medskip
 The property of quasi-periodicity of \ $\zeta_4$, \ i.e.
formula (5) of \SS 4, allows us to obtain other formulas.

\hskip 0,45cm First,  let us introduce :
$$\eqalign{\eta (x,\omega )
&= 2 \zeta_4 (\omega ) + \bigl( (x+\omega ) \mid \nabla_\omega )^2
\zeta_4 (\omega )\cr
&= 2\bigl( 1 + ( (x+\omega ) \mid \nabla_\omega )^2\bigr) \zeta_4
(\omega ),\cr}$$

so that the quasi-periodicity of \ $\zeta_4$ \ could be written as :
$$\zeta_4 (x+2\omega ) = \zeta_4 (x) + \eta (x, \omega )\leqno (6)$$

\bigskip
Look now at some elementary properties of \ $\eta (x, \omega )$~: 
\medskip
{\parindent=1cm

\item{(i)} $\eta$ is a polynomial on \ $x$ \ of degree 2 over \ $\R_{0,3}$.
\item{} Take \ $\zeta_4(x)$.  We can write on one hand :
$$\zeta_4 (x) = \zeta_4 (x-2\omega + 2\omega ) = \zeta_4
(x-2\omega ) + \eta (x, \omega )$$
\item{} so that \ $\zeta_4 (x-2\omega ) = \zeta_4(x) - \eta (x, \omega )$. 
But, on the other hand
$$\zeta_4 (x-2\omega ) = \zeta_4 (x) + \eta (x, - \omega ).$$

\bigskip\bigskip
\item{} In such a way we get :
\smallskip
\item{(ii)} $\eta (x, -\omega ) = - \eta (x, \omega )$,  which means that \
$\eta$ \ is a even function on its second variable.
\item{} The argument below gives also :
\par}
$$\zeta_4 (x-2\omega ) = \zeta_4 (x) - \eta (x, \omega ).\leqno (7)$$ 

{\parindent=1cm
Now compute :
$$\eqalign{\eta (-x, \omega )
&= 2\zeta (\omega ) + \bigl( (-x+\omega ) \mid \nabla_\omega )^2 \zeta
(\omega ) \cr
&= -2\zeta (-\omega ) - \bigl( (x-\omega ) \mid \nabla_\omega )^2
\zeta (-\omega )\cr
&= - \eta (x, -\omega ) = \eta (x, \omega ) \cr}$$

\bigskip\medskip
which shows :
\smallskip
\item{(iii)} $\eta$  is an even function of its first variable $x$.

\bigskip\medskip
\item{} As a direct consequence of (ii) and (iii) we get :
\smallskip
\decale{(iv)} $\eta (-x, - \omega ) = - \eta (x, \omega )$.

\item{}The last result could be obtained also by a direct computation
carried on the defining expression of $\eta$.
\item{}Developping the defining expression for \ $\eta$ \ we have~:
$$\eqalign{\eta (x, \omega )
&= 2\zeta (\omega ) + (\omega \mid \nabla_\omega )^2 \zeta (\omega )\cr
&+ 2 (x \mid \nabla_\omega ) (\omega \mid \nabla_\omega ) \zeta
(\omega )\cr
&+ (x \mid \nabla_\omega )^2 \zeta (\omega ),\cr}$$

\bigskip\medskip
\item{}in which we have to erase those terms which are zero, i.e. those in
which appears \ $(\omega\mid\nabla_\omega ) \zeta (\omega )$ \ so~:
\medskip
\item{(v)} $\eta (x, \omega ) = 2\zeta (\omega ) +
(x\mid\nabla_\omega )^2 \zeta (\omega )$.

\bigskip\medskip
In the last, set \ $x = \omega$.  Thus we have :
\item{(vi)} $\eta (\omega , \omega ) = 2\zeta (\omega )$.
\par}

\bigskip
Remark \ (vi) \ could be obtained directly from (6) setting \ $x = -\omega$,
\ or from (7) with \ $x = \omega$.

\bigskip
We will end this paragraph with the following : take two half-periods \
$\omega^1$ \ and \ $\omega^2$.  Set \ $\omega^{12}$ \ such that
$$\omega^1 + \omega^2 + \omega^{12} = 0.$$

Apply formula (7) on :
$$\eqalign{
&\zeta_4 (x+2\omega^{12}) = \zeta_4 (x-2\omega^1 - 2\omega^2) =\cr
&= \zeta_4 (x-2\omega^1) - \eta (x, \omega^2) =\cr
&= \zeta_4 (x) - \eta (x, \omega^1) - \eta (x, \omega^2).\cr}$$

On the other hand, from (6) we have :
$$\zeta_4(x+2\omega^{12}) = \zeta_4(x) + \eta (x, \omega^{12}),$$
so that :
$$\eta (x, \omega^1) + \eta (x, \omega^2) + \eta  (x, \omega^{12}) = 0.$$

\vskip 1cm
\font\grandsy=cmsy10 scaled\magstep2
\def\SSg{{\grandsy x}}
{\bf \SSg\kern.15em  6. Cliffordian analogues of the Weierstrass \
${\cal P}$ \ function}
\bigskip\medskip
\font\grandsy=cmsy10 scaled\magstep0
\def\SS{{\grandsy x}}
 Remember the classical situation of elliptic functions. Once we
have constructed the \ $\zeta$ \ function, which is an odd function,
quasi-periodic on $2\omega_1$ \ and \ $2\omega_2$,  and which is defined as~:
$$\zeta (z) = {1\over z} + \sum_{p=1}^\infty \ \left\{ {1\over z-w_p} +
{1\over w_p} + {z\over x_p^2}\right\} ,$$

where the set \ $(2\Z\omega_1 + 2\Z\omega_2) \setminus \{ (0,0)\}$ \ has
been rearranged in a set denoted by \ $\{ w_p\}$, $p\in\N^*$, then, in order
to introduce an authentic 2-periodic function, we need to derive $\zeta$~:
$$\zeta'(z) = {d\over dz} \zeta (z) = {\partial\over\partial x} \zeta
(z) = - i {\partial\over\partial y} \zeta (z).$$

The tradition demands to set :
$${\cal P} (z) = - \zeta' (z),$$
so that :
$${\cal P} (z) = {1\over z^2} + \sum_{p=1}^\infty \ \left\{ {1\over
(z-w_p)^2} - {1\over w_p^2}\right\} .$$

\bigskip
Obviously, \ ${\cal P}$ \ is an even function. One want to show that \
${\cal P}$ \ is a 2-periodic function. For this purpose there is a
traditionnal argument based on the fact that obviously \ ${\cal P}'$ \ is a
2-periodic function and because of the parity of \ ${\cal P}$ \ one can
conclude, by an integration, that \ ${\cal P}$ \ is also 2-periodic.

\bigskip
 In our situation, we start with the \ $\zeta_4$ \ function~:
$$\zeta_4(x) = x^{-1} + \sum_{p=1}^\infty \{ (x-w_p)^{-1} + \sum_{\mu = 0}^3
(w_p^{-1}x)^\mu \ w_p^{-1} \} .$$

Here, the possibilities of derivations are larger than in the previous case,
so the candidates for analogues of the \ ${\cal P}$ \ function are numerous.

\bigskip\medskip
 It would be possible to restrict us on a canonical way of
derivation, just on the real axis, namely on \
$\displaystyle{\partial\over\partial x_0}$.  Taking into account we need
three successive derivations, one of the \ ${\cal P}_4$ \ functions would
be~:
$$\leqalignno{ {\cal P}_4(x)
&= - {1\over 3!} \ {\partial^3\over \partial x_0^{3} }\  \zeta _4(x) =&(8)\cr
&= x^{-4} + \sum_{p=1}^\infty \{ (x-w_p)^{-4} - w_p^{-4}\} .&\cr}$$

\bigskip
In fact, the right notation for the last function should be \ ${\cal
P}_{(3,0,0,0)} (x)$, \ in which we indicate the used derivations.

\bigskip\bigskip
 Generally, we have 20 \ ${\cal P}$-functions coming from \
$(e_i \mid \nabla_x)^3 \ \zeta_4(x)$, $i = 0,1,2,3$.  In the same way
we introduced the notation for the holomorphic Cliffordian polynomials in
\SS 1, if \ $\alpha = (\alpha_0, \alpha_1, \alpha_2, \alpha_3)$ \ is a
multiindex in  \ $\N^4$ \ with lenght \ $\vmid{\alpha} = 3$, \ we can look
at the \ ${\cal P}_\alpha (x)$ \ functions.

\bigskip
 For example \  :

\hskip 1cm ${\cal P}_{(3,0,0,0)} (x) = (e_0 \mid\nabla_x)^3 \ 
\zeta_4(x)$, 

\hskip 1cm ${\cal P}_{(2,0,1,0)} (x) = (e_0\mid\nabla_x)^2 (e_2
\mid\nabla_x) \  \zeta_4(x)$,

\hskip 1cm ${\cal P}_{(0,1,1,1)}(x) = (e_1\mid
\nabla_x) (e_2 \mid \nabla_x) (e_3 \mid \nabla_x) \ 
\zeta_4(x)$. 

\bigskip
 Now let us show that the \ ${\cal P}_\alpha$ \ functions are 
4-periodics. For this we will observe that, again by derivation, one
gets~:
$$(e_i \mid \nabla_x) {\cal P}_\alpha (x+2\omega ) = (e_i \mid
\nabla_x) {\cal P}_\alpha (x).$$

Here \ $i = 0,1,2,3$, \ $\alpha = (\alpha_0, \alpha_1, \alpha_2, \alpha_3)$,
\ $\vmid{\alpha } \  = 3$.

\bigskip\medskip
\hskip 0,45cm Integrating this four relations respectively in the directions
\ $e_0, e_1, e_2, e_3$, \ one deduces~:
\medskip
$${\cal P}_\alpha (x+2\omega ) = {\cal P}_\alpha (x) + C_i (x_0,\ldots ,
\hat x_i,\ldots, x_3)$$

for \ $i = 0,1,2,3$, \ where \ $C_i (x_0,\ldots,\hat x_i,\ldots,x_3)$ \
designe \ ``constants" \ depending on the mentioned variables. By successive
substractions we can prove the independance of the variables of the
constants. So for every \ $x\in S\oplus V$, \ one has~: 
$${\cal P}_\alpha (x+2\omega ) = {\cal P}_\alpha (x) + C$$
with a constant $C$. Setting \ $x = -\omega$, we have :
$${\cal P}_\alpha (\omega ) = {\cal P}_\alpha (-\omega ) + C$$
from which we deduce \ $C = 0$ \ because of the parity of \ ${\cal
P}_\alpha$.  Thus the 4-periodicity of \ $P_\alpha$ \ was proved.

\bigskip
 All the \ ${\cal P}_\alpha$ \ functions are elliptic
Cliffordian with pointwise poles of order 4 at the vertices of the latice \
$2\Z^4\omega$.

\vskip 1cm
\font\grandsy=cmsy10 scaled\magstep2
\def\SSg{{\grandsy x}}
{\bf \SSg\kern.15em  7. On the sets of the zeros and the poles of \
$D_0 {\cal P}_0$ }
\bigskip\medskip
 Here \ $D_0 = \displaystyle{\partial\over\partial x_0}$ \ and \
${\cal P}_0$ \ is an abreviated notation for \ ${\cal P}_{(3,0,0,0)}$.
\medskip

 According to  (8) it is clear that the set of the poles of
\ $D_0 {\cal P}_0$ \ is \ $2\Z^4\omega$.

\bigskip
 Let now restrict us in an elementary cell in \ $\R^4$ \ which
is nothing else than a hypercube spanned over four linearly independant
paravectors \ $2\omega_1, 2\omega_2, 2\omega_3, 2\omega_4$ \ (in fact four \
$\R$-independant vectors in \ $\R^4$).  Thus, the only pole of \ $D_0P_0$ \
in this cell is a pole of order 5 at the origin.

\bigskip
 Try now to find the zeros of \ $D_0 {\cal P}_0$.  Remember \
$D_0 {\cal P}_0$ \ is 4-periodic and odd.  So~:
$$D_0 {\cal P}_0 (x+2\omega ) = D_0 {\cal P}_0(x).$$

Set \ $x = -\omega$,  one deduces
$$D_0 {\cal P}_0 (\omega ) = D_0 {\cal P}_0 (-\omega ) = - D_0 {\cal P}_0
(\omega )$$

from which :
$$D_0 {\cal P}_0 (\omega ) = 0$$

that means \ $D_0 {\cal P}_0$ \ vanishes at all the half-periods \
$\omega_1, \omega_2, \omega_3, \omega_4$ \ and their linear combinations by
2, by 3 elements and finally at the vertex \
$\omega_1+\omega_2+\omega_3+\omega_4$ \ that is at all the vertices
(excepting $0$)  of the hypercube spanned over \ $\omega_1, \omega_2,
\omega_3, \omega_4$. The number of this vertices is \ $2^4-1 = 15$.

\bigskip
 At this stage we do not have guarantees this are the only zeros
of \ $D_0 {\cal P}_0$. The result is only that, concerning the 4-periodic
function \ $D_0 {\cal P}_0$,
$${\cal Z} \geq 3P$$
where \ ${\cal Z}$ \ means the set of its zeros and \ $P$ \ the set of his
poles, as usually,  taking into account the multiplicities.

\bigskip\medskip
\undertext{Conjecture}.  Is is true that in an elementary cell of the
4-periodic lattice, every 4-periodic function \ $f$ \ satisfies
$${\cal Z}_f = 3 P_f~?$$

\vskip1,5cm
\centerline{\bf Bibliographie}
\bigskip\medskip

\bb 1&F. BRACKS, R. DELANGHE, F. SOMMEN&Clifford analysis,&Pitman, (1982) & &

\bb 2&C.A. DEAVORS&The quaternion calculus&Am. Math. Monthly. (1973), 995-1008& &

\bb 3&R. DELANGHE, F. SOMMEN, V. SOUC\v EK&Clifford Algebra and
Spinor-valued functions&Kluwer Academic Publishers& &

\bb4&R. FUETER&Die Funktionnentheorie der Differentialgleichungen \ $\Delta
u = 0$ und $\Delta\Delta u = 0$ mit vier reellen Variablen.& Comment Math. Helv 7
(1935), 307-330& &

\bb5&R. FUETER&Uber die analytische Darstellung der regularen
Funktionen einer Quaternionenvariabelen&Comm. Math. Helv.8 (1936), 371-378& &

\bb6&G. LAVILLE&Une famille de solutions de l'\'equation de Dirac avec champ
\'electroma\-gn\'etique quelconque&C.R. Acad. Sci. Paris t. 296 (1983),
1029-1032& &

\bb7&G. LAVILLE&Sur l'\'equation de Dirac avec champ \'electromagn\'etique
quelconque&Lectures Notes in Math. 1165, Springer-Verlag (1985), 130-149& &

\bb8&G. LAVILLE, I. RAMADANOFF&Fonctions holomorphes Cliffordiennes&C.R. Acad,
Sc. Paris, 326, s\'erie I (1998), 307-310& &

\bb9&G. LAVILLE, I. RAMADANOFF&Holomorphic Cliffordian Functions&Advances in
Applied Clifford Algebras, 8, n$^\circ$2 (1998), 323-340& &

\bb10&H. MALONEK&Powers series representation for monogenic functions in \
$\R^{n+1}$ based on a permutational product&Complex variables, vol 15 (1990),
181-191& &

\bb11&L. PERNAS&Holomorphie quaternionienne&Advances in Applied Clifford
Algebras, 8, n$^\circ$2 (1998), 283-298& &

\bb12&P.L. WALKER&Elliptic Functions&J. Wiley and  Sons (1996)& &

\bigskip

{\parindent=8cm
\item{UPRES-A 6081} D\'epartement de Math\'ematiques
\item{}Universit\'e de Caen
\item{}14032 CAEN Cedex France

\medskip
\item{}glaville@math.unicaen.fr
\item{}rama@math.unicaen.fr
\par}

\end